\documentclass[english,11pt]{article}
\usepackage[latin1]{inputenc} 
\usepackage{amsmath,amsthm,amssymb,babel}
\renewcommand{\arraystretch}{1.7}
\begin{document}
\newtheorem{teo}{Theorem}
\newtheorem{defin}{Definition}
\newtheorem{prop}{Proposition}
\newtheorem{cor}{Corollary}
\newtheorem{lemma}{Lemma}
\newtheorem{rem}{Remark}
\renewcommand{\arraystretch}{1.7}
\def\proof{{\it Proof.}\ }
\def\endproof{\hfill $\Box$\par\vskip3mm}
\def\pt#1{{\it Proof of Theorem \ref{#1}.}}
\def\eq#1{(\ref{#1})}
\def\th#1{Theorem \ref{#1}}
\def\neweq#1{\begin{equation}\label{#1}}
\def\endeq{\end{equation}}
\def\weak{\rightharpoonup}
\def\ep{\varepsilon}
\def\half{{1\over2}}
\def\phi{\varphi}
 \def\RR{{\mathbb R} }
\def\CC{{\mathbb C} }
\def\NN{{\mathbb N} }
\def\supp{{\rm supp}}
\def\di{\displaystyle}
\def\ri{\rightarrow}
\def\intom{\int_\Omega}

\title{\sc Entire solutions of sublinear elliptic equations in anisotropic media}
\author{Teodora-Liliana Dinu \\ \small Department of Mathematics, ``Fra\c tii Buze\c sti" College, 
200352 Craiova, Romania\\ \small E-mail: {\tt tldinu@gmail.com} }

\date{}
\maketitle

\noindent{\small{\sc Abstract}.
We study the nonlinear elliptic problem $-\Delta u=\rho (x)f(u)$ in $\RR^N$ ($N\geq 3$),
$\lim_{|x|\ri\infty}u(x)=\ell$,
where $\ell\geq 0$ is a real number, $\rho(x)$ is a nonnegative potential belonging to a certain 
Kato class,
and $f(u)$ has a sublinear growth.
We distinguish the cases $\ell>0$ and $\ell=0$ and we prove existence and uniqueness results if 
the potential $\rho(x)$ decays fast enough at infinity. Our arguments rely on comparison 
techniques and on a theorem of Brezis and Oswald for sublinear elliptic equations.  
}

\section{Introduction and main results}
In their celebrated paper \cite{bk}, Brezis and Kamin have been concerned with various questions 
related to the existence of bounded solutions of the sublinear elliptic equation without condition 
at infinity
\begin{equation}\label{brek}
-\Delta u=\rho (x)u^\alpha\qquad\mbox{in $\RR^N$, $N\geq 3$},\end{equation}
where $0<\alpha<1$, $\rho\in L^\infty_{\rm loc}(\RR^N)$, $\rho\geq 0$ and $\rho\not\equiv 0$.
We summarize in what follows the main results obtained in \cite{bk}.
Brezis and Kamin proved that the {\it nonlinear} problem \eq{brek}  has a bounded solution $u> 0$ 
if and only if the {\it linear} problem
$$-\Delta u=\rho (x)\qquad\mbox{in $\RR^N$}$$
has a bounded solution. In this case, Problem \eq{brek} has a minimal positive solution
and this solution satisfies $\liminf_{|x|\ri\infty}u(x)=0$. Moreover, the minimal solution is the 
unique positive solution of \eq{brek} which tends to zero at infinity. Brezis and Kamin also 
showed that if the potential $\rho(x)$ decays fast enough at infinity then Problem \eq{brek} has a 
solution and, moreover, such a solution does not exist if $\rho (x)$ has a slow decay at infinity. 
For instance, if $\rho(x)=(1+|x|^p)^{-1}$, then \eq{brek} has a bounded solution if and only if 
$p>2$.
More generally, Brezis and Kamin have proved that Problem \eq{brek} 
has a bounded
solution if and only if $\rho(x)$ is potentially bounded, that is, the mapping
$x\longmapsto \int_{\RR^N}\rho(y)|x-y|^{2-N}dy\in L^\infty(\RR^N)$. We
refer to \cite{1,kra} for various results on bounded domains for sublinear elliptic equations with 
zero Dirichlet boundary condition. Problem \eq{brek} in the whole space has been considered in 
\cite{eu,ede,egn,2,kav,mab,nai} under various assumptions on $\rho$. Sublinear problems (either 
stationary or evolution ones) appear in the study of population dynamics, of reaction-diffusion 
processes, of filtration in porous media with absorption, as well as in the study of the scalar 
curvature of warped products of semi-Riemannian manifolds (see, e.g., \cite{on}).

Our purpose in this paper is to study the problem
\begin{equation}
\label{ec}
\left\{
\begin{tabular}{ll}
$-\Delta u=\rho(x)f(u)$ &in $\RR^N$\\
$u>\ell$ &in $\RR^N$\\
$u(x)\ri \ell$ &as $|x|\ri\infty$, 
\end{tabular}\right.
\end{equation}
where $N\geq 3$ and $\ell\geq 0$ is a real number. 

Throughout the paper we assume that the variable potential $\rho(x)$ satisfies $\rho\in 
L^\infty_{\rm loc}(\RR^N)$, $\rho\geq 0$ and $\rho\not\equiv 0$. 

In our first result we suppose
that the growth at infinity of the anisotropic potential $\rho(x)$  is given by

\smallskip
($\rho$1) $\int _0^{\infty}r\Phi(r)dr<\infty$,
where $\Phi (r):=\max _{|x|=r}\rho(x).$

\smallskip
Assumption ($\rho$1) has been first introduced in Naito \cite{nai}.

The nonlinearity $f:(0,\infty)\ri (0,\infty)$ satisfies $f\in C^{0,\alpha}_{\rm loc}(0,\infty)$
($0<\alpha<1$)
and has a sublinear growth, in the sense that 

\smallskip
(f1) the mapping $u \longmapsto f(u)/u$ is decreasing on
$(0,\infty )$ and $\lim_{u\rightarrow \infty}f(u)/u
=0$.

\smallskip
We point out that condition (f1) does not require that $f$ is smooth at the origin. The standard 
example of such a nonlinearity is $f(u)=u^p$, where $-\infty<p<1$. We also observe that we study 
an equation of the same type as in Brezis and Kamin \cite{bk}. The main difference is that we 
require a certain asymptotic behaviour at infinity of the solution.


Under the above hypotheses ($\rho$1) and (f1), 
our first result concerns the case $\ell>0$. We have

\begin{teo}
\label{t1} Assume that $\ell>0$.
Then Problem \eq{ec} has a unique classical solution.
\end{teo}

Next, consider the case $\ell=0$. Instead of ($\rho$1) we impose the stronger condition

\smallskip
($\rho$2) $\int _0^\infty r^{N-1}\Phi (r)dr <\infty$.

\smallskip
We remark that in Edelson \cite{ede2} it is used the stronger assumption
$\int _0^\infty r^{N-1+\lambda(N-2)}\Phi (r)dr <\infty$, for some $\lambda\in (0,1)$. 

Additionally, we suppose that

\smallskip
(f2) $f$ is increasing in $(0,\infty)$ and $\lim_{u\searrow 0}f(u)/u=+\infty$.

\smallskip
A nonlinearity satisfying both (f1) and (f2) is $f(u)=u^p$, where $0<p<1$.

Our result in the case $\ell=0$ is the following.

\begin{teo}\label{t2} Assume that $\ell=0$ and assumptions ($\rho$2),
(f1) and (f2) are fulfilled. Then Problem \eq{ec} has
a unique classical solution.
\end{teo}

We point out that assumptions ($\rho$1) and ($\rho$2) are related to a 
celebrated class introduced by
Kato, with wide and deep applications in Potential Theory and Brownian
Motion. We recall (see \cite{aize}) that a real-valued measurable
function $\psi$ on $\RR^N$ belongs to the Kato class ${\cal K}$
provided that
$$\lim_{\alpha\ri
0}\sup_{x\in\RR^N}\int_{|x-y|\leq\alpha}E(y)|\psi (y)|dy=0,$$
where $E$ denotes the fundamental solution of the Laplace
equation. According to this definition and our assumption ($\rho$1) (resp., ($\rho$2)),
it follows that $\psi=\psi(|x|)\in {\cal K},$ where
$\psi (|x|):=|x|^{N-3}\Phi(|x|)$ (resp., $\psi (|x|):=|x|^{-1}\Phi(|x|)$), for all $x\not=0.$

\section{Proof of Theorem \ref{t1}}

In order to prove the existence of a solution to Problem \eq{ec}, we use
 a result established by
 Brezis and  Oswald (see \cite[Theorem~1]{1})
 for bounded domains. Consider the problem
\begin{equation}\label{1}
\left\{\begin{array}{lll}
\displaystyle -\Delta u=g(x,u) \quad\mbox{in\ }
\Omega \\
\displaystyle u\geq 0,\ \ \ \,u \not\equiv 0
\quad\mbox{in\ } \Omega \\
\displaystyle \, \, \ \, u=0 \quad\mbox{on\ }
\partial\Omega \, ,
\end{array}\right.
\end{equation}
where $ \Omega\subset \RR^{N} $ is a bounded domain
with smooth boundary and
$ g(x,u):\Omega \times [0,\infty) \rightarrow \RR$.
 Assume that
\begin{equation}\label{2}
 \left\{\begin{array}{lll}
 \displaystyle  \mbox{\rm for\ a.e.\ } x\in \Omega \
\mbox{\ the\ function\ }u\longmapsto g(x,u) \
\mbox{ is continuous on }\ [0,\infty) \\
 \displaystyle \mbox {\rm\ and\ the\ mapping\ } u\longmapsto g(x,u)/u
 \mbox{\rm\ is\ decreasing\ on\ }(0,\infty)\, ;
\end{array}  \right.
\end{equation}
\begin{equation}\label{3}
\begin{array}{lll}
\mbox{\ for\ each\ } u\geq 0 \mbox{\ the\ function\ }
x\longmapsto
g(x,u) \mbox{\ belongs\ to\ } L^{\infty}(\Omega)\, ;
\end{array}
\end{equation}
\begin{equation}\label{4}
  \exists C>0 \mbox{\ such that\ }
 g(x,u)\leq C(u+1) \mbox{\ a.e.}\ x\in \Omega,\ \ \ \forall\ u\geq 0.
\end{equation}
Set
$$ a_{0}(x)=\lim\limits_{u\searrow 0}g(x,u)/u \qquad \mbox{and}
\qquad
a_{\infty}(x)=\lim\limits_{u\ri \infty}g(x,u)/u\, ,  $$
so that $ -\infty<a_{0}(x)\leq +\infty $ and
$ -\infty\leq a_{\infty}(x)<+\infty. $

Under these hypotheses, Brezis and Oswald proved
in \cite{1} that Problem \eq{1}
has at most one solution.
   Moreover, a solution
   of \eq{1} exists if and only if
   \begin{equation}\label{108}
\lambda_{1}(-\Delta-a_{0}(x))<0
   \end{equation}
and
   \begin{equation}\label{109}
\lambda_{1}(-\Delta-a_{\infty}(x))>0,
   \end{equation}
where $ \lambda_{1}(-\Delta-a(x)) $
denotes the first eigenvalue of the operator
$ -\Delta-a(x) $
with zero Dirichlet condition.
The precise meaning of $ \lambda_{1}(-\Delta-a(x)) $
 is
$$ \lambda_{1}(-\Delta-a(x)) = \inf_{\varphi\in
H_{0}^{1}(\Omega),\ \|\varphi\|_{L^2(\Omega)}=1}\left( \int |\nabla
\varphi|^{2}-\int\limits_{[\varphi\not=0]}a
\varphi^{2} \right). $$
Note that $ \displaystyle
\int\limits_{[\varphi\not=0]} a \varphi^{2} $
makes sense if $ a(x) $ is any measurable
function such that either $ a(x)\leq C $ or $ a(x)\geq -C $ a.e. on $ \Omega. $

For any positive integer $k$ we consider the problem
\begin{equation}
\label {margptc}
\left\{
\begin{tabular}{ll}
$-\Delta u_k=\rho(x)f(u_k)$ , &{\rm if} $|x|<k$\\
$u_k>\ell$, &{\rm if} $|x|<k$\\
$u_k(x)=\ell,$           &{\rm if } $|x|=k.$
\end{tabular}\right.
\end{equation}
Equivalently, the above boundary value problem can be rewritten 
\begin{equation}
\label {marg}
\left\{
\begin{tabular}{ll}
$-\Delta v_k=\rho(x)f(v_k+\ell)$ , &{\rm if} $|x|<k$\\
$v_k(x)=0,$           &{\rm if } $|x|=k.$
\end{tabular}\right.
\end{equation}
In order to obtain a solution of the problem \eq{marg}, it is enough to check the
hypotheses of the Brezis- Oswald theorem.
\begin{itemize}
\item Since $f\in C(0,\infty)$ and $\ell>0$, it follows that the mapping $v\longmapsto 
\rho(x)f(v+\ell)$ is continuous in $[0,\infty)$.
\item
From $\rho(x)\frac {f(v+\ell)}{v}=\rho(x)\frac{f(v+\ell)}{v+\ell}\frac{v+\ell}{v}$,
using positivity of $\rho$ and (f1) we deduce that the function $v\longmapsto
\rho(x)\frac {v+\ell}{v}$ is decreasing on $(0,\infty)$.
\item
For all $v\geq 0$, since $\rho\in L^{\infty}_{\rm loc}(\RR^N)$, we obtain that $\rho\in L^{\infty}
(B(0,k))$, so the condition \eq{3} is satisfied.
\item
By $\lim_ {v\rightarrow \infty}f(v+\ell)/(v+1)=0$ and $f\in C(0,\infty)$,
there exists $M>0$ such that $f(v+\ell)\leq M(v+1)$ for all $v\geq 0.$ Therefore
$\rho(x)f(v+\ell)\leq ||\rho||_{L^{\infty}(B(0,k))}M(v+1)$ for all $v\geq 0.$
\item We have
$$a_0(x)=\lim_ {v\searrow 0}\frac {\rho(x)f(v+\ell)}{v}=+\infty$$
and
$$a_{\infty}(x)=\lim_ {v\rightarrow \infty}\frac{\rho(x)f(v+\ell)}{v}=
\lim_ {v\rightarrow \infty}\rho(x)\frac{f(v+\ell)}{v+\ell}\cdot\frac{v+\ell}{v}=0\,.$$
\end{itemize}
Thus, by Theorem 1 in \cite{1},
Problem \eq{marg} has a unique solution $v_k$ which, by the
maximum principle, is positive in $|x|<k.$  Then $u_k=v_k+\ell$ satisfies \eq{margptc}.
Define $u_k=\ell$ for $|x|>k.$ The maximum principle
implies that
$\ell\leq u_k\leq u_{k+1}$ in $\RR^N.$

We now justify the existence of a continuous function $v:\RR^N\ri\RR$, $v>\ell$, such that 
 $u_k\leq v$ in $\RR^N.$ We first construct a
positive radially symmetric function $w$ such that $-\Delta
w=\Phi(r)$ $(r=|x|)$ in $\RR^N$ and $\lim _{r\rightarrow
\infty}w(r)=0.$ A straightforward computation shows that
$$w(r)=K-\int _0^r\zeta ^{1-N}\int _0^\zeta \sigma ^{N-1}\Phi(\sigma)d\sigma d\zeta,$$
where $$K=\int _0^{\infty}\zeta ^{1-N}\int _0^\zeta \sigma
^{N-1}\Phi(\sigma)d\sigma d\zeta,$$ provided the integral is
finite. An integration by parts yields
$$\begin{array}{ll}\di\int _0^{r}\zeta ^{1-N}\int _0^\zeta \sigma ^{N-1}\Phi(\sigma)d\sigma 
d\zeta=&\di
-\frac{1}{N-2}\int _0^r \frac{d}{d\zeta}\zeta ^{2-N}\int_0^\zeta\sigma ^{N-1}\Phi(\sigma)
d\sigma d\zeta\\
&\di=\frac{1}{N-2}\left(-r^{2-N}\int_0^r \sigma^{N-1}\Phi(\sigma)d\sigma+
\int_0^r\zeta \Phi(\zeta)\right)\\
&\di <\frac{1}{N-2}\int_0^{\infty}\zeta \Phi(\zeta)<+\infty\,.\end{array}$$
Moreover, $w$ is decreasing and satisfies $0<w(r)<K$ for all $r\geq 0$.
Let $v>\ell$ be a function such that $w(r)=m^{-1} \int _0^{v(r)-\ell}\frac{t}{f(t+\ell)}dt,$
where $m>0$ is chosen such that $Km\leq \int_0^m \frac{t}{f(t+\ell)}dt.$

Next, by L'H\^opital's rule for the case $\frac{.}{\infty}$ (see
\cite[Theorem 3, p.~319]{mn}) we have
$$\lim \limits _ {x\rightarrow \infty}\di\frac{\int _0^x\frac{t}{f(t+\ell)}dt}{x}=
\lim \limits _ {x\rightarrow \infty}\frac{x}{f(x+\ell)}=
\lim \limits _ {x\rightarrow \infty}\frac{x+\ell}{f(x+\ell)}\cdot \frac{x}{x+\ell}=+\infty\,.$$
This means that there exists $x_1>0$ such that $\int _0^x\frac{t}{f(t)}\geq Kx$
for all $x\geq x_1$. It follows that for any $m\geq x_1$ we have $Km\leq \int _0^m
\frac{t}{f(t)}dt.$

Since $w$ is decreasing, we obtain that $v$ is a decreasing function, too. Then
$$\int _0^{v(r)-\ell}\frac{t}{f(t+\ell)}dt\leq \int _0^{v(0)-\ell}\frac{t}{f(t+\ell)}dt=mw(0)=
mK\leq \int _0^m \frac{t}{f(t+\ell)}dt\,.$$
It follows that $v(r)\leq m+\ell$ for all $r>0.$

From $w(r)\rightarrow 0$ as $r\rightarrow \infty$ we deduce that $v(r)\rightarrow \ell$
as $r\rightarrow \infty.$

By the choice of $v$ we have
$$\nabla w=\frac{1}{M}\frac{v-\ell}{f(v)}\nabla v \quad {\rm and}\quad \Delta w=
\frac {1}{m}\frac{v-\ell}{f(v)}\Delta v+\frac{1}{m}\left(\frac {v-\ell}{f(v)}\right)'|\nabla 
v|^2\,.$$
Since the mapping $u\longmapsto f(u)/u$ is decreasing  on $(0,\infty)$
we deduce that
\begin{equation}
\label {relv}
\Delta v<\frac{m}{v-\ell}f(v)\Delta w=-\frac{M}{v-\ell}f(v)\Phi(r)\leq -f(v)\Phi(r)\,.
\end{equation}
By \eq{margptc}, \eq{relv} and our
hypothesis (f1), we
obtain that $u_k(x)\leq v(x)$ for each $|x|\leq k$ and so, for all $x\in\RR^N.$

In conclusion,
$$u_1\leq u_2\leq \ldots\leq u_k\leq u_{k+1}\leq \ldots\leq v,$$
with $v(x)\rightarrow \ell$ as $|x|\ri\infty$. Thus, there exists a function
 $u\leq v$ such that $u_k\rightarrow u$ pointwise in $\RR^N.$ In particular, this shows that
 $u>\ell$ in $\RR^N$ and $u(x)\rightarrow \ell$ as $|x|\ri\infty$.

A standard bootstrap argument (see, e.g., \cite{gt}) shows that  $u$ is a classical solution of 
the problem \eq{ec}. 

To conclude the proof, it remains to show that the solution found above is unique.
Suppose that $u$ and $v$ are solutions of \eq{ec}.
It is enough to show that $u\leq v$ or, equivalently, $\ln u(x)\leq \ln
v(x),$ for any $x\in \RR^N.$ Arguing by contradiction, there exists
$\overline x\in \RR$ such that $u(\overline x) >v(\overline x).$
Since $\lim_ {|x|\rightarrow \infty}(\ln u(x)-\ln v(x))=0$,
we deduce that $\max _{\RR^N} (\ln u(x)-\ln v(x))$ exists and is positive. At this point,
say $x_0$, we have
\begin{equation}
\label{nab}
\nabla (\ln u(x_0)-\ln v(x_0))=0,
\end{equation}
so
\begin{equation}
\label{ineg}
\frac{\nabla u(x_0)}{u(x_0)}=\frac{\nabla v(x_0)}{v(x_0)}.
\end{equation}
By (f1) we obtain
$$\frac {f(u(x_0))}{u(x_0)}<\frac {f(v(x_0))}{v(x_0)}.$$
So, by \eq{nab} and \eq{ineg},
$$\begin{array}{ll}\di 0&\di\geq \Delta (\ln u(x_0)-\ln v(x_0))\\
&\di =\frac{1}{u(x_0)}\cdot \Delta u(x_0)-
\frac{1}{v(x_0)}\cdot \Delta v(x_0)
-\frac {1}{u^2(x_0)}\cdot |\nabla u(x_0)|^2+\frac {1}{v^2(x_0)}\cdot |\nabla v(x_0)|^2\\
&\di =
\frac{\Delta u(x_0)}{u(x_0)}-\frac {\Delta v(x_0)}{v(x_0)}=
-\rho(x_0)\left(\frac{f(u(x_0))}{u(x_0)}-\frac{f(v(x_0))}{v(x_0)}\right)>0,\end{array}$$
which is a contradiction. Hence $u\leq v$ and the proof is concluded. \qed

\section{Proof of Theorem \ref{t2}}
\subsection{Existence}
Since $f$ is an increasing positive function on $(0,\infty)$, there exists
and is finite $\lim_ {u\searrow 0}f(u),$ so $f$ can be extended by continuity
at the origin.
Consider the Dirichlet problem
\begin{equation}
\label {marg0}
\left\{
\begin{tabular}{ll}
$-\Delta u_k=\rho(x)f(u_k)$ , &{\rm if} $|x|<k$\\
$u_k(x)=0,$           &{\rm if } $|x|=k.$
\end{tabular}\right.
\end{equation}
Using the same arguments as in case $\ell>0$ we deduce that conditions \eq{2} and \eq{3}
are satisfied. In what concerns assumption \eq{4}, we use both assumptions (f1) and (f2). Hence
$f(u)\leq f(1)$ if $u\leq 1$ and $f(u)/u\leq f(1)$ if $u\geq 1.$ Therefore $f(u)\leq f(1)(u+1)$ 
for all $u\geq 0$, which proves \eq{4} .
The existence of a solution for \eq{marg0} follows from \eq{108} and \eq{109}. These conditions 
are direct consequences of our assumptions
$\lim_{u\rightarrow \infty}f(u)/u=0$ and
$\lim_ {u\searrow 0}f(u)/u=+\infty$. Thus, by the Brezis-Oswald theorem, Problem
\eq{marg0} has a unique solution.
Define $u_k(x)=0$ for $|x|>k.$
Using the same arguments as in case $\ell>0,$
we obtain $ u_k\leq u_{k+1} $ in $ \RR^N.$

Next, we prove the existence of a continuous function
$ v:\RR^N\ri\RR $ such that $ u_k\leq v $ in $ \RR^N.$
We first construct a positive radially
symmetric function $w$
satisfying
 $ -\Delta w=\Phi(r) $ ($r=| x|$) in $ \RR^N $
and
 $  \lim_{r\to \infty}w(r)=0. $
We obtain
$$ w(r)=\displaystyle K - \int\limits_0^r \zeta^{1-N}\int
\limits_0^{\zeta}\sigma^{N-1}\Phi(\sigma) d\sigma d\zeta \, ,$$
where
\begin{equation}\label{kmare}
 K=\displaystyle
\int\limits_0^{\infty}\zeta^{1-N}\int\limits_0^{\zeta}
\sigma^{N-1}\Phi(\sigma)\,d\sigma\,d\zeta   \, ,
\end{equation}
provided the  integral is finite.
By integration by parts we have 
 \begin{equation}\label{9}
   \begin{array}{lll}
\di\int\limits_0^r \zeta^{1-N}\int\limits_0^{\zeta}\sigma^{N-1}
\Phi(\sigma)\,d\sigma\,d\zeta =-\frac{1}{N-2}\int\limits_0^r
\frac{d}{d\zeta}\zeta^{2-N}\int\limits_0^{\zeta}
\sigma^{N-1}\Phi(\sigma)\,d\sigma\,d\zeta  = \\
\di \frac{1}{N-2}\left(-r^{2-N}\int\limits_0^r\sigma^{N-1}
 \Phi(\sigma)\,d\sigma+\int\limits_0^r \zeta  \Phi(\zeta)\,d\zeta \right) <
\frac{1}{N-2} \int_0^{\infty}\zeta  \Phi(\zeta)\,d\zeta<\infty\,.
 \end{array}    \end{equation}
Therefore
$$ w(r)<\frac{1}{N-2}\cdot \int\limits_0^{\infty} \zeta \Phi(\zeta)d\zeta,\qquad\mbox{
for all $r>0$} .$$
Let $ v $ be a positive function such that $ w(r)=
c^{-1}\int_0^{v(r)}t/f(t)dt, $ where $ c>0 $
is chosen
such that $ K c\leq  \int_{0}^{c}t/f(t)dt.$
 We argue in what follows that we can find $ c>0 $ with this property. Indeed, by
 L'H\^opital's rule,
$$ \lim_{x\to \infty}\frac{\int\limits_{0}^{x}
\frac{t}{f(t)}dt}{x}=\lim_{x\to \infty}\frac{x}{f(x)}=
+\infty\, . $$
This means that there exists $ x_1>0 $ such that
$ \int_{0}^{x}t/f(t)dt\geq Kx $
for all $ x\geq x_{1}. $ It follows that for any $ c\geq x_1 $ we
have $ K c\leq \int_{0}^{c}t/f(t)dt.$

On the other hand, since $w$ is decreasing, we deduce that $v$ is a decreasing function, too.
Hence $$ \int\limits_{0}^{v(r)}\frac{t}{f(t)}dt\leq \int\limits_{0}^{v(0)}\frac{t}{f(t)}dt=c\cdot 
w(0)=c\cdot K\leq \displaystyle \int\limits_{0}^{c}\frac{t}{f(t)}dt\,.$$
It follows that $ v(r)\leq c $ for all $ r>0 $.

 From $ w(r)\to 0 $ as $ r\to \infty $ we deduce that
 $ v(r)\to 0$ as $r\to \infty .$

 By the choice of $ v$ we have
    \begin{equation}\label{10}
\nabla w=\frac{1}{c}\cdot \frac{v}{f(v)}\nabla v \
\mbox{\ and\ }\
\Delta w=\frac{1}{c} \frac{v}{f(v)}\Delta v +
\frac{1}{c}\left(\frac{v}{f(v)}\right)'|\nabla v|^2 .
    \end{equation}
Combining the fact that 
$ f(u)/u $ is a
decreasing function on $ (0,\infty)$ with relation \eq{10}, we deduce that
\begin{equation}\label{nr}
 \Delta v < c\frac{f(v)}{v}\Delta w=-c \,
\frac{f(v)}{v}\Phi(r)\leq -f(v)\Phi(r)\, .
\end{equation}
By \eq{marg0} and \eq{nr} and using  our hypothesis
(f2), as already done for proving the uniqueness in the case $\ell>0$, we
obtain that $ u_k(x)\leq v(x) $ for each $ |x|\leq k $ and so,
for all $x\in \RR^{N}. $

 We have obtained a bounded increasing sequence
$$ u_1\leq u_2\leq \ldots \leq u_k\leq u_{k+1}\leq
 \ldots \leq v \, ,$$
with $ v$ vanishing at infinity. Thus, there exists a
function $ u\leq v$ such that $ u_k\to u $
pointwise in $ \RR^N. $
A standard bootstrap argument implies that
 $u$ is a classical solution of the problem \eq{ec}.

\subsection{Uniqueness}
We split the proof into two steps. Assume that $u_1$ and $u_2$ are solutions of Problem
\eq{ec}. We first prove that if $u_1\leq u_2$ then $u_1=
u_2$ in $\RR^N.$ In the second step we find a positive
solution $u\leq \min\{u_1,u_2\}$ and thus, using the first
step, we deduce that  $u= u_1$ and $u= u_2$, which
proves
the uniqueness.

\medskip
{\sc Step I}. We show that $u_1\leq u_2$ in $\RR^N$ implies $u_1= u_2$ in
$\RR^N.$ Indeed, since
$$u_1\Delta u_2-u_2\Delta u_1=\rho(x)u_1u_2\left(\frac{f(u_2)}{u_2}-\frac {f(u_1)}{u_1}\right)
\geq 0\,,$$
it is sufficient to check that
\begin{equation}
\label {intpeRn}
\int _{\RR^N}(u_1\Delta u_2-u_2\Delta u_1)=0
\end{equation}
Let $\psi \in C_0^{\infty}(\RR^N)$ be such that
$\psi(x)=1$ for $|x|\leq 1$ and $\psi(x)=0$ for $|x|\geq 2$, and
denote $\psi _n:=\psi (x/n)$ for any positive integer $n$.
Set
$$I_n:=\int _{\RR^N}(u_1\Delta u_2-u_2\Delta u_1)\psi _ndx.$$
We claim that $I_n\rightarrow 0$ as $n\rightarrow \infty$. Indeed,
$$|I_n|\leq \int _{\RR^N}|u_1\Delta u_2|\psi _ndx+\int _{\RR^N}|u_2\Delta u_1|\psi _ndx\,.$$
So, by symmetry, it is enough to prove that $J_n:=\int
_{\RR^N}|u_1\Delta u_2|\psi _ndx\rightarrow 0$ as $n
\rightarrow \infty.$  But, from \eq{ec},
\begin{equation}\label{ineg2}
\begin{array}{ll}\di J_n&\di=\int _{\RR^N}|u_1f(u_2)\rho(x)|\psi _ndx=\int _n^{2n}\int \limits 
_{|x|=r}
|u_1(x)f(u_2(x))\rho(x)|dxdr\\ &\di\leq 
\int _n^{2n}\Phi(r) \int \limits_{|x|=r}|u_1(x)f(u_2(x))|dxdr\leq
\int _n^{2n}\Phi(r) \int \limits_{|x|=r}|u_1(x)|M(u_2+1)dxdr\,.\end{array}
\end{equation}
Since $u_1(x)$, $u_2(x)\rightarrow 0$ as $|x|\rightarrow \infty$, we deduce
that $u_1$ and $u_2$ are bounded in $\RR^N$. Returning to \eq{ineg2} we have
$$J_n\leq M(||u_2||_{L^\infty(\RR^N)}+1)\sup _{|x|\geq n}|u_1(x)|\cdot \frac{\omega_N}{N}
\int _n^{2n}\Phi(r)r^{N-1}dr
\leq C\int _0^{\infty}\Phi (r)r^{N-1}dr\cdot \sup _{|x|\geq n}|u_1(x)|\,.$$
Since $u_1(x)\rightarrow 0$ as $|x|\rightarrow \infty$, we have
$\sup_{|x|\geq n}|u_1(x)|\rightarrow 0$ as $n\rightarrow
\infty$ which shows that $J_n\rightarrow 0.$ In particular, this implies
$I_n \rightarrow 0$ as $n\ri\infty$. 

We recall in what follows the Lebesgue Dominated Convergence Theorem (see \cite[Theorem 
IV.2]{brezis}).
\begin{teo}
\label{leb} Let $f_n:\RR^N\ri\RR$ be a sequence of  functions in $L^1(\RR^N)$. We assume that

(i)  $f_n(x)\ri f(x)$ a.e. in $\RR^N$,

(ii) there exists $g\in L^1(\RR^N)$ such that, for all $n\geq 1$, 
$|f_n(x)|\leq g(x)$ a.e. in $\RR^N$.

 Then $f\in L^1(\RR^N)$ and $||f_n-f||_{L^1}\ri 0$ as $n\ri\infty$.
\end{teo}

Taking $f_n:=(u_1\Delta u_2-u_2\Delta u_1)\psi _n$
 we deduce $f_n(x)\rightarrow u_1(x)\Delta
u_2(x)-u_2(x)\Delta u_1(x)$ as $n\rightarrow \infty$. To apply
Theorem \ref{leb} we need to show that $u_1\Delta u_2-u_2\Delta
u_1 \in L^1(\RR^N).$ For this purpose  it is
sufficient to prove that $u_1\Delta u_2\in L^1(\RR^N).$ Indeed,
$$\int _{\RR^N}|u_1\Delta u_2|\leq ||u_1||_{L^\infty(\RR^N)}\int _{\RR^N}|\Delta u_2|=C
\int _{\RR^N}|\rho(x)f(u_2)|\,.$$ Thus, using $f(u)\leq f(1)(u+1)$ and since $u_2$ 
is bounded,
the above inequality yields
$$\begin{array}{ll}\di\int _{\RR^N}|u_1\Delta u_2|&\di\leq C \int _{\RR^N}|\rho(x)(u_2+1)|\\
&\di\leq 
C\int _0^{\infty}\int _{|x|=r}\Phi(r)dxdr\leq C \int 
_0^{\infty}\Phi(r)r^{N-1}<+\infty\,.\end{array}$$
This shows that $u_1\Delta u_2\in L^1(\RR^N)$ and the proof of Step I is completed.

\medskip
{\sc Step II}. Let $u_1$, $u_2$ be arbitrary solutions of Problem \eq{ec}. For all integer $k\geq 
1$, denote $\Omega_k:=\{x\in\RR^N;\ |x|< k\}$. The Brezis-Oswald theorem implies that the problem
$$
\left\{
\begin{tabular}{ll}
$-\Delta v_k=\rho(x)f(v_k)$ &{\rm{in}} $\Omega _k$\\
$v_k=0$ &{\rm{on}} $\partial \Omega _k$
\end{tabular}\right.
$$
has a unique solution $v_k\geq 0.$ Moreover, by the
Maximum Principle, $v_k>0$ in
$\Omega_k.$ We define $v_k=0$ for $|x|>k.$ Applying again the Maximum
Principle we deduce that $v_k\leq v_{k+1}$ in $\RR^N.$
Now we prove that $v_k\leq u_1$ in $\RR^N$, for all $k\geq 1.$
Obviously, this happens  outside $\Omega_k.$ On the other hand
$$
\left\{
\begin{tabular}{ll}
$-\Delta u_1=\rho(x)f(u_1)$ &{\rm{in}} $\Omega_k$\\
$u_1>0$ &{\rm{on}} $\partial \Omega_k$
\end{tabular}\right.
$$
Arguing by contradiction, we assume that there exists $\overline x\in\Omega_k$ such that 
$v_k(\overline{x})>u_1(\overline x).$
Consider the function $h:\Omega_k \rightarrow \RR$, $h(x)=\ln v_k(x)-\ln u(x).$
Since $u_1$ is bounded in $\Omega_k$ and $\inf _{\partial \Omega _k}u_1>0$ we have
$\lim \limits _ {|x|\rightarrow k}h(x)=-\infty$. We deduce that
$\max _{\Omega_k}\left(\ln v_k(x)-\ln u_1(x)\right)$ exists and is positive. Using the same 
argument
as in the case $\ell>0$ we deduce that $v_k\leq u_1$ in $\Omega _k$, so in $\RR ^N.$
Similarly we obtain $v_k\leq u_2$ in $\RR ^N.$ Hence
$v_k\leq \overline u:=\min\{u_1,u_2\}$.
Therefore $v_k\leq v_{k+1}\leq \ldots\leq \overline u$.
 Thus there exists a
function $u$ such that $v_k\rightarrow u$ pointwise in
$\RR^N.$ Repeating a previous argument  we
deduce that $u\leq \overline u$ is a classical solution of Problem
\eq{ec}. Moreover, since $u\geq
v_k>0$ in
$\Omega _k$ and for all $k\geq 1$, we deduce that $u>0$ in $\RR^N$. This concludes the proof of 
Step II.

Combining Steps I and II we conclude that  $u_1= u_2$ in
$\RR^N.$ \qed

\end{document}